\documentclass{article}%
\usepackage{amsmath}%
\usepackage{amsfonts}%
\usepackage{amssymb}%

\newtheorem{thm}{Theorem}
\newtheorem{pro}[thm]{Proposition}

\newtheorem{rem}{Remark}[section]

\def\l{\lambda}
\def\s{\sigma}
\def\proof{{\noindent \bf Proof. \hspace{.01in}}}
\newcommand{\eop}{\hspace{.1in} \vrule height 7pt width 5pt depth 0pt \medskip}

\begin{document}


\title{\Large Harnack type inequalities for matrices in majorization}

\author{Chaojun Yang${}^{\,\rm a}$
,\;
Fuzhen Zhang${}^{\,\rm b}$
\\
\footnotesize{Department of Mathematics}\\
\footnotesize{${}^{\,\rm a}$ Suzhou University, Suzhou, China; cjyangmath@163.com}\\
\footnotesize{${}^{\,\rm b}$ Nova Southeastern University, Fort Lauderdale, USA; zhang@nova.edu}}

\date{}
\maketitle

\hrule
\bigskip


\noindent {\bf Abstract.}
Following the recent work of Jiang and Lin (Linear Algebra Appl. 585 (2020) 45--49), we present more results (bounds) on Harnack type  inequalities  for matrices in terms of majorization
(i.e., in partial products) of eigenvalues and singular values. We   discuss and compare  the bounds derived through different ways.   Jiang and Lin's results imply Tung's version of Harnack's inequality
(Proc. Amer. Math. Soc. 15 (1964)  375--381); our results 
are stronger and more general than
Jiang and Lin's. We also show some majorization inequalities 
concerning Cayley transforms. Some open problems on spectral norm and eigenvalues are proposed.
\medskip

{\footnotesize

}

\bigskip
\hrule


\section{Introduction}
There are several mathematical inequalities that carry Harnack's name in the literature.
The classical Harnack inequality is about relating the values of a positive harmonic function at two points in a domain. The inequality is usually shown
by using Poisson's formula with integration on a sphere; see \cite{Kass} for a nice introduction about the inequality and its proof.
Generalized Harnack inequalities in various forms have been developed and heavily used in  partial differential equations \cite{Di12,  MorTian07, Mueller06, WFY}.

We are  concerned with the Harnack inequality for matrices.

Tung \cite{Tun64}  established the following determinantal Harnack inequality.

\begin{thm}[Tung]\label{Tung}
 Let $Z$ be an $n\times n$ complex matrix with singular values $r_k$  that satisfy  $0\le r_k<1$,   $k=1, 2, \ldots, n$ {\rm (}i.e., $Z$ is a strict contraction{\rm )}. Let $Z^*$ denote the conjugate transpose of $Z$ and let $I$ be the $n\times n$ identity matrix.  Then for any $n\times n$ unitary matrix  $U$, it holds true that
\begin{eqnarray}\label{tung1}
\prod_{k=1}^n\frac{1-r_k}{1+r_k}\le \frac{\det(I-Z^*Z)}{|\det(I-UZ)|^2}\le \prod_{k=1}^n\frac{1+r_k}{1-r_k}.
\end{eqnarray}
Equality occurs on the right if and only if
$UZ$ has eigenvalues  $r_1, r_2, \dots, r_n$;   equality on the left holds if and only if
$UZ$ has eigenvalues  $-r_1, -r_2, \dots, -r_n$ \nolinebreak {\rm \cite{LinZ17}}.
\end{thm}

Proved by a Language multiplier method, Tung's work drew immediate  attention of
Hua   and Marcus.  Hua \cite{Hua65} gave a proof of    (\ref{tung1}) using a determinantal inequality
 he had previously obtained in \cite{Hua55}, while Marcus \cite{Mar65} considered an equivalent form
   of  (\ref{tung1})  without denominators.
   Later, Fan \cite{KyFan78, FanLAA1988} formulated and proved Harnack's inequalities  for operators with norm less than 1
    in the setting of  Hilbert space.
   Recent work on the matrix Harnack inequality includes \cite{LinZ17} in which the inequality is extended to multiple contractive matrices  and \cite{JiangLAA20} in which the authors present more general forms of (\ref{tung1}) in majorization sense \cite{MOA11}. For recent work on contractive matrices, see
\cite{Omar19}.

With $A=UZ$, (\ref{tung1}) is  equivalently rewritten  in terms of eigenvalues as
\begin{eqnarray}\label{tung2}
\prod_{k=1}^n\frac{1-r_k}{1+r_k}\le\prod_{k=1}^n  \l_k\Big (
(I-A^*)^{-1}(I-A^*A) (I-A)^{-1}\Big ) \le \prod_{k=1}^n\frac{1+r_k}{1-r_k}.
\end{eqnarray}

The matrix in the product in the middle of  (\ref{tung2})  automatically gets attention as it is
a term in the Schur complement of
$\left ( {(I-A^*A)^{-1} \atop (I-A^*)^{-1} }
{(I-A)^{-1} \atop \cdot } \right )$  which  resembles the Hua matrix
$\left ( {(I-A^*A)^{-1} \atop (I-A^*B)^{-1} }
{(I-B^*A)^{-1} \atop \cdot } \right )$   
  \cite{And80, FanLAA1988, XXZLAMA11, XXZLAA09}  and as the Julia operator
  $\left ( {(I-AA^*)^{\frac12} \atop -A^* } {A \atop (I-A^*A)^{\frac12} }
 \right )$ is unitary  \cite[p.\,148]{You88}.
The latter two block matrices (operators) have often been used in deriving matrix or operator inequalities.

(\ref{tung2})  leads to the study of inequalities of  partial products, i.e., log-majorization,
of eigenvalues and singular values. (Note that inequalities in log-majorization are  in general stronger than (weak-) majorization inequalities which are equivalent to the inequalities in
unitarily invariant norms.)
Following this line, an interesting generalization of (\ref{tung1}) is presented by Jiang and Lin in \cite{JiangLAA20}. Our goal   is to  continue with Jiang and Lin's work and to show more results (bounds) of this type. We   compare  the bounds derived through  different approaches.

\section{Main results}

We state our first result for matrices. The identities in fact  hold true
for  linear operators in a complex Hilbert space. Let $M_n$ be the space of $n\times n$ complex matrices. For $X\in M_n$, let
$\Re (X)=\frac12 (X+X^*)$ and $\Im (X)=\frac{1}{2i} (X-X^*)$, where $X^*$ is the adjoint (conjugate transpose) of $X$. $X=\Re (X)+i \Im (X)$ is the Cartesian decomposition of $X$.
Let $\Lambda(X)$ denote the spectrum of $X$.

 \begin{pro}\label{Thm:exp1}
 Let $A\in M_n$
  such that $1\not \in \Lambda(A)$. Then
\begin{eqnarray}
\lefteqn{(I-A^*)^{-1}(I-A^*A)(I-A)^{-1}} \nonumber \qquad   \\
  & = & 2\, \Re \big ((I-A)^{-1}\big ) -I \label{expz} \\
 & = & 2 \,\Re \big ((I-A)^{-1}-{\small{\frac{1}{2}} }I\big  ) \label{exp2} \\
  & = & \Re \big ( (I+A)(I-A)^{-1}\big  ) \label{fan}\\
  & = & S^*S,\, S=(I-A^*A)^{\frac12}(I-A)^{-1} \mbox{if $A$ is contractive}.\label{exp3}  
  \end{eqnarray}
  \end{pro}

  \proof The first identity, i.e., (\ref{expz}), is the same as
  $$(I-A^*)^{-1}(I-A^*A)(I-A)^{-1}=(I-A)^{-1}+(I-A^*)^{-1}-I$$
  which is easily verified by  multiplying by $I-A^*$ from the left and by $I-A$ from the right.  (\ref{exp2}) is immediate from (\ref{expz}).   (\ref{fan}) holds true if and only if
 $$2(I-A^*)^{-1}(I-A^*A)(I-A)^{-1}=(I+A)(I-A)^{-1}+(I-A^*)^{-1}(I+A^*), $$
equivalently, by multiplying by $I-A^*$ from the left and by $I-A$ from the right,
 $$2(I-A^*A)=(I-A^*)(I+A)+(I+A^*)(I-A)$$
 which is obvious. (\ref{exp3}) is trivial. \eop

 The identity or expression  (\ref{fan}) in Proposition \ref{Thm:exp1}
  appeared in    \cite{FanLAA1988};  it was used as a pivot in \cite{JiangLAA20} to obtain  the desired inequalities.
  Fan derived the identity using analysis with assumption $\|A\|<1$. This condition is unnecessary in (\ref{expz})--(\ref{fan}). We will obtain various bounds by the expressions in Proposition~\ref{Thm:exp1}.

  For $X\in M_n$,  let $\Lambda (X)=\{\l_1(X), \dots, \l_n(X)\}$ be the set of the  eigenvalues of $X$. The eigenvalues are arranged in non-increasing order if they are all real, i.e.,  $\l_1(X)\geq \l_2(X)\geq  \cdots \geq \l_n(X)$. For singular values,   we denote by  $\s_j(X)$ the $j$th largest singular value of $X$, i.e.,
$\s_j(X)=\sqrt{\l_j(X^*X)}$,
and $\s_1(X)\geq \s_2(X)\geq  \cdots \geq \s_n(X)$.
 For simplicity, sometimes we use $r_j$ for $\s_j(\cdot )$.

Our main theorem is on the upper and lower bounds with singular values.
\begin{thm}\label{thm:main}
 Let $A\in M_n$ be a strict contraction {\rm (}which implies  $1\not \in \Lambda(A)${\rm )} with singular values ordered as
$0\leq r_n\leq \cdots \leq r_2\leq r_1<1$. Then  for   $k=1, 2, \dots, n$, and for any sequence
$1\leq i_1< \cdots < i_k\leq n$, the following inequalities hold:
\begin{equation}\label{J0}
\l_{j} \big ( (I-A^*)^{-1}(I-A^*A)(I-A)^{-1} \big )  \leq \frac{1+r_{j}}{1-r_{j}}, \;\;
 j=1, 2, \dots, n,
\end{equation}
\begin{equation}\label{J1}
\prod_{j=1}^k \l_{i_j} \Big ( (I-A^*)^{-1}(I-A^*A)(I-A)^{-1} \Big ) \leq \prod_{j=1}^k\frac{1+r_{i_j}}{1-r_{i_j}}\leq \prod_{j=1}^k\frac{1+r_{j}}{1-r_{j}},
\end{equation}
\begin{equation}\label{Jb}
\prod_{j=1}^k \l_{n-{i_j}+1} \Big ( (I-A^*)^{-1}(I-A^*A)(I-A)^{-1} \Big ) \geq \prod_{j=1}^k\frac{1-r^2_{i_j}}{(1+r_{j})^2}\geq \prod_{j=1}^k\frac{1-r_{j}}{1+r_{j}}.
\end{equation}
\end{thm}

\proof  To prove (\ref{J0}), we borrow two  known facts: For any   $X\in M_n$,
\smallskip

(i). $\l_j (\Re (X))\leq \s_j(X)$, $j=1, 2, \dots, n$, and

(ii). $\s_j(X)+\s_{n-j+1}(I-X)\geq 1$; $\s_i(X)+\s_{j}(I-X)\geq 1$ if $i+j\leq n+1$.
\smallskip

Fact (i) is a well-known result of Fan and Hoffman \cite{FanHoffman55}, while (ii) is immediate from  the fact that $\s_{i+j-1}(X+Y)\leq \s_i(X)+\s_j(Y)$ for $i+j\leq n+1$ applied to  $I=X+(I-X)$. See, e.g., \cite[pp.\,73--75]{BhaMA97}. It is obvious that in (ii)  we can replace $I$ by any $n\times n$  unitary matrix. In addition, $\s_j(I-X)\leq 1+\s_j(X)$.

We now use   (i), (ii), and expression (\ref{expz}) in Proposition \ref{Thm:exp1} to derive
\begin{eqnarray*}
\mbox{The left-hand side (LHS) of (\ref{J0})} & = &  \l_{j} \big (2\, \Re ((I-A)^{-1})-I\big )\\
& = &  2\,\l_{j}  \big  (  \Re ((I-A)^{-1})\big )-1\\
 & \leq &  2\,\s_{j}     ((I-A)^{-1})-1 \;\; 
 \\
 & = & \frac{2}{\s_{n-{j}  +1} (I-A)}-1\\
  & \leq  & \frac{2}{1-\s_{{j}  } (A)}-1\\
  & = &
   \frac{2}{1-r_j} -1 \\
  & =& \frac{1+r_{j}}{1-r_{j}}.
 \end{eqnarray*}

Notice that every term in (\ref{J0}) is positive when $A$ is a strict contraction. Thus, (\ref{J1}) follows immediately from (\ref{J0}).
  The   last inequality in (\ref{J1})  is due to the fact that $f(t)=\frac{1+t}{1-t}$ is an increasing function on $[0, 1).$
\medskip

 To prove (\ref{Jb}), we use the following results: 
 For any $A, B\in M_n$,
\smallskip

(iii). $\prod_{j=1}^k\s_{i_j}(AB)\geq  \prod_{j=1}^k\s_{n-j+1}(A)\s_{i_j}(B)$, $k=1, 2, \dots, n$,
and

(iv). $\prod_{j=1}^k\s_{i_j}(AB)\leq  \prod_{j=1}^k\s_j(A)\s_{i_j}(B)$, $k=1, 2, \dots, n$.

\smallskip

(iv) is well-known, see, e.g., \cite[p.\,72]{BhaMA97}, \cite[p.\,340]{MOA11}, or \cite[p.\,364]{ZFZbook11}. (iii) and (iv) are in fact equivalent. (iv)
 implies (iii)  by replacing
$A$ with $A^{-1}$ and $B$ with $AB$. (If $A$ is singular, then use a continuity argument.)

\smallskip
Now we use  (\ref{exp3}) in Proposition  \ref{Thm:exp1} and  compute
\begin{eqnarray*}
\mbox{LHS  of (\ref{Jb})} & = &
\prod_{j=1}^k \l_{n-{i_j}+1}  (S^*S )= \prod_{j=1}^k \s^2_{n-{i_j}+1}(S) \\
& = & \Big (\prod_{j=1}^k \s_{n-{i_j}+1}\big ( (I-A^*A)^{\frac12}(I-A)^{-1}\big )\Big )^2\\
 & \geq & \Big (\prod_{j=1}^k \s_{n-{i_j}+1}\big ( (I-A^*A)^{\frac12}\big ) \, \s_{n-j+1}\big ((I-A)^{-1}\big )\Big )^2
 \\
 & = &  \prod_{j=1}^k \frac{1-\s^2_{i_j}(A)}{\s^2_{j} (I-A)}
   \geq   \prod_{j=1}^k \frac{1-\s^2_{i_j}(A)}{\big (1+\s_{j} (A)\big )^2}\\
    &  = &  \prod_{j=1}^k\frac{1-r^2_{i_j}}{(1+r_{j})^2}\geq   \prod_{j=1}^k\frac{1-r^2_{j}}{(1+r_{j})^2}
   =  \prod_{j=1}^k\frac{1-r_{j}}{1+r_{j}}.  \eop
 \end{eqnarray*}

Setting $i_j=j$ in (\ref{J1}) and (\ref{Jb}) reveals the inequalities     in \cite{JiangLAA20}.

 \begin{rem} \rm On upper bounds.
By  (\ref{exp3}) and (iv), we get another upper bound:
\begin{equation}\label{J2}
\prod_{j=1}^k \l_{i_j} \Big ( (I-A^*)^{-1}(I-A^*A)(I-A)^{-1} \Big ) \leq
\prod_{j=1}^k
 \frac{1-r_{n-j+1}^2}{(1-r_{{i_j}})^2}.
\end{equation}

The proof goes as follows.
\begin{eqnarray*}
\mbox{LHS of (\ref{J2})}& = & \prod_{j=1}^k \s_{i_j}^2 \big ((I-A^*A)^{\frac12}(I-A)^{-1}\big )\\
& \leq  & \Big (  \prod_{j=1}^k \s_j \big ((I-A^*A)^{\frac12}\big )\, \s_{i_j}  \big ((I-A)^{-1}\big )  \Big )^2\\
 & = &  \Big (  \prod_{j=1}^k \big (1-\s_{n-j+1}^2(A)\big )^{\frac12} \cdot  \frac{1}{\s_{n-{i_j}+1} (I-A)}  \Big )^2 \\
  & \leq  & \prod_{j=1}^k
 \frac{1-\s_{n-j+1}^2(A)}{\big (1-\s_{{i_j}} (A)\big )^2}\\
 & = & \prod_{j=1}^k
 \frac{1-r_{n-j+1}^2}{(1-r_{{i_j}})^2}.
 \end{eqnarray*}

In a similar way, using (iv), one   obtains  the  upper bound
 $\prod_{j=1}^k
 \frac{1-r_{n-_{i_j}+1}^2}{(1-r_{j})^2}$ in place of
 $\prod_{j=1}^k \frac{1-r_{n-j+1}^2}{(1-r_{{i_j}})^2}$ in (\ref{J2}).
Comparisons of these  bounds 
 are in order.
Let
{\small
$$
 R_1=\prod_{j=1}^k\frac{1+r_{i_j}}{1-r_{i_j}}, \; R_2=\prod_{j=1}^k\frac{1+r_{j}}{1-r_{j}}, \; R_3=\prod_{j=1}^k
 \frac{1-r_{n-j+1}^2}{(1-r_{{i_j}})^2}, \; R_4=\prod_{j=1}^k
 \frac{1-r_{n-_{i_j}+1}^2}{(1-r_{j})^2}.$$}

 We saw $R_1\leq R_2$ in the proof of Theorem~\ref{thm:main}. We claim $R_1\leq R_3$, $R_2\leq R_4$, but $R_2$ and $R_3$ are incomparable, and $R_3$ and $R_4$ are incomparable.

Let $a_i=1-r_{n-j+1}^2$. Then $a_1\geq a_2\geq \cdots \geq a_n\geq 0$, and the product $a_1a_2\cdots a_k$ is greater than or equal to the product of any $k$ of $a$'s.
 It follows that
 $$R_1=\prod_{j=1}^k\frac{1-r_{i_j}^2}{(1-r_{i_j})^2}\leq \prod_{j=1}^k
 \frac{1-r_{n-j+1}^2}{(1-r_{{i_j}})^2}=R_3.$$
 For a similar reason by considering the product of $k$ smallest $a$'s, we get $R_2\leq R_4$.

If $r=(r_1, r_2, r_3, r_4)=(\frac12, \frac12, 0, 0)$, $k=1$, $i_1=3$, then $R_2=3 > 1=R_3$;
if  $r=(\frac12, \frac{9}{20}, 0, 0)$, $k=1$, $i_1=2$, then $R_2=3 < \frac{400}{121} =R_3<4=R_4$.
To have $R_4<R_3$, we take  $r=(\frac12, \frac12, \frac12, 0, 0)$, $k=2$, $i_1=2,$  $ i_2=3$. Then
$R_3=16>12=R_4$. Thus, $R_2$ and $R_3$ are incomparable, so are $R_3$ and $R_4$.

Therefore,
$R_1\leq R_2\leq R_4, \, R_1\leq R_3.$  Moreover, if  we set  $R_5=\prod_{j=1}^k
 \frac{1-r_{n-j+1}^2}{(1-r_{{j}})^2}$, it is easy to show that
$R_3\leq R_5$ and $R_4\leq R_5$. Of all the upper bounds obtained above  in  log-majorization, we conclude that $R_1$ is optimal.
\end{rem}

\begin{rem}\rm On lower bounds. As
$\frac{1+r_j}{1-r_j}$ is an upper bound in (\ref{J0}),
it is natural and interesting  to ask if the reversal $\frac{1-r_j}{1+r_j}$ can serve as a lower bound.

From the proof of (\ref{J0}), we see
the upper bound  essentially follows from  the inequality $\l_j(\Re((I-A)^{-1}))\leq \frac{1}{1-r_j}$.
It is tempting to have $\frac{1}{1+r_j}$ as a lower bound for $\l_j(\Re((I-A)^{-1}))$ that  would result in the lower bound $\frac{2}{1+r_j}-1=\frac{1-r_j}{1+r_j}$ in~(\ref{J0}).
However, this is not true in general.
Take
$$A=\left ( \begin{array}{ccc}
 0.4831    & 0.2041   &  0.0447\\
    0.4689 &   0.3308 &   0.3671\\
    0.1308  &  0.2583  &  0.4787 \end{array} \right ).$$
    Then the singular values of $A$ are  $
      0.9468, $ $
    0.3969, $ $
    0.0049$,
  the eigenvalues of $\Re ((I-A)^{-1})$ are $
    9.9860, $ $
    1.5616, $ $
    0.7789$,
    and the eigenvalues of
$(I-A^*)^{-1}(I-A^*A)(I-A)^{-1}$ are
   $ 18.9720, $ $
    2.1232, $ $
    0.5578.$

   One may check that $\l_3(\Re((I-A)^{-1}))=0.7789<0.9951=\frac{1}{1+r_3}$, and
   $$\l_3\big ((I-A^*)^{-1}(I-A^*A)(I-A)^{-1}\big )=0.5578<0.9902=\frac{1-r_3}{1+r_3}.$$

 Moreover,
$\prod_{j=1}^k\frac{1-r^2_{i_j}}{(1+r_{j})^2}$ in (\ref{Jb})
can be similarly  replaced by
$\prod_{j=1}^k\frac{1-r^2_{j}}{(1+r_{i_j})^2}$. Setting $k=1$ and replacing $i_j$ by $n-j+1$,  we arrive at, for   $j=1, 2, \dots, n$,
\begin{equation}\label{JbX}
\l_{j} \Big ( (I-A^*)^{-1}(I-A^*A)(I-A)^{-1} \Big ) \geq
\frac{1-r^2_{n-j+1}}{(1+r_{1})^2}\geq \frac{1-r_{1}}{1+r_{1}}.
\end{equation}
and
\begin{equation}\label{JbY}
\l_{j} \Big ( (I-A^*)^{-1}(I-A^*A)(I-A)^{-1} \Big ) \geq
\frac{1-r^2_{1}}{(1+r_{n-j+1})^2}\geq \frac{1-r_{1}}{1+r_{1}}.
\end{equation}

The previous example shows that  the middle terms in (\ref{JbX}) and (\ref{JbY}) cannot be replaced by $\frac{1-r_j}{1+r_j}$. But can the $r_1$'s in (\ref{JbX}) and (\ref{JbY}) be replaced by $r_{n-j+1}$?  See the later (\ref{open2}) and  the $j$-conjecture in the next section.
 \end{rem}

\section{Fan's norm inequalities and open problems} 
Let $\|A\|$ denote the spectral (operator) norm of a bounded linear operator  $A$ on a complex Hilbert space.
For  $A$ with $\|A\|<1$, Fan \cite[Prop.\,1 (3)]{FanLAA1988} showed that
\begin{equation}\label{FanProp1(3)}
\frac{1-\|A\|}{1+\|A\|}(I-A^*)(I-A)
\leq I-A^*A\leq
\frac{1+\|A\|}{1-\|A\|}(I-A^*)(I-A),
\end{equation}
where $H\leq K$ means that $H$, $K$ are self-adjoint and $K-H$ is a positive operator.
(Note that the above inequalities (\ref{FanProp1(3)}) (i.e., (3) in \cite{FanLAA1988}) imply other inequalities in  Proposition 1 of Fan \cite{FanLAA1988}. For instance, one can derive the left inequality of (1) of  Fan \cite{FanLAA1988} from
the left inequality of  (\ref{FanProp1(3)}); and vice versa, as Fan showed.)

In case of matrices,  $\|A\|$ is equal to   the largest singular value of $A$, i.e., $r_1$ in the previous sections.  It follows that  Fan's (\ref{FanProp1(3)}) is equivalent to
\begin{equation}\label{FanProp1(3b)}
\frac{1-r_1}{1+r_1} I
\leq (I-A^*)^{-1}(I-A^*A)(I-A)^{-1}\leq
\frac{1+r_1}{1-r_1} I.
\end{equation}

(\ref{FanProp1(3b)})
follows from   Theorem \ref{thm:main} (with $k=1$) immediately because
$$\frac{1-r_1}{1+r_1}\leq \l_{j} \big ( (I-A^*)^{-1}(I-A^*A)(I-A)^{-1} \big )   \leq
\frac{1+r_1}{1-r_1}.$$

Theorem \ref{thm:main} (\ref{J0}) presents    stronger   upper bounds,  for   $j=1, 2, \dots, n$,
$$
\l_{j} \big ( (I-A^*)^{-1}(I-A^*A)(I-A)^{-1} \big )  \leq \frac{1+r_j}{1-r_j}\leq
\frac{1+r_1}{1-r_1}.$$

For   $A\in M_n$, let $|A|=(A^*A)^{1/2}$. Observing  that
$$\frac{1+r_j}{1-r_j}=\frac{2}{1-r_j}-1=\l_j\big( 2(I-|A|)^{-1}-I\big ),$$  we can rewrite (\ref{J0}) as, for   $j=1, 2, \dots, n$,
$$
\l_{j} \big ( (I-A^*)^{-1}(I-A^*A)(I-A)^{-1} \big )  \leq \l_j\big( 2(I-|A|)^{-1}-I\big ),$$
or equivalently,
$$
\l_{j} \big ( (I-A^*)^{-1}+(I-A)^{-1}  \big ) \leq  \l_j\big(  2(I-|A|)^{-1}\big ),$$
i.e.,
$$
\l_{j} \big ( \Re ((I-A)^{-1})  \big ) \leq  \l_j\big( (I-|A|)^{-1}\big ).$$
It is natural to ask if the stronger inequalities in the Loewner sense
hold:
$$
 (I-A^*)^{-1}(I-A^*A)(I-A)^{-1}  \leq  2(I-|A|)^{-1}-I,$$
 or equivalently,
 $$
 (I-A^*)^{-1}+(I-A)^{-1}  \leq  2(I-|A|)^{-1}.$$
 This is false in general as one may verify  with
 $A=\left ( {0 \atop 0}{0.1\atop 0}\right )$  that
$$ (I-A^*)^{-1}+(I-A)^{-1}=\left ( \begin{array}{cc}
2 &  0.1 \\
  0.1 &   2
   \end{array} \right )\not \leq 2(I-|A|)^{-1}=\left ( \begin{array}{cc}
2 &  0  \\
  0 &   \frac{20}{9}
   \end{array} \right ).$$
The same example also shows that
$$ (I-A^*)^{-1}+(I-A)^{-1} \not \geq 2(I+|A|)^{-1}.$$
   However, a great amount of numerical computation shows  that for each $j$,
\begin{equation}\label{open}
\l_j\big ((I-A^*)^{-1}+(I-A)^{-1}\big )    \geq 2\l_j\big ((I+|A|)^{-1}\big )=
\frac{2}{1+r_{n-j+1}},
\end{equation}
or equivalently,
\begin{equation}\label{open2}
\l_j\big ((I-A^*)^{-1}(I-A^*A)(I-A)^{-1}\big )    \geq \l_j\big (2 (I+|A|)^{-1}-I\big )=\frac{1-r_{n-j+1}}{1+r_{n-j+1}}.
\end{equation}

If (\ref{open}) and (\ref{open2})  hold true, then we would have nice  lower bounds for (\ref{J0}).

 We propose two open problems; the second one is a special case of the first.
  Let $A$ be an $n\times n$ strict contraction, i.e., the spectral  norm $\|A\|<1$. Then
\begin{equation}\label{jConj}
\l_j\big (\Re((I-A)^{-1})\big )    \geq \l_j\big ((I+|A|)^{-1}\big ),  \;\; j=1, 2,  \dots, n.
\end{equation}

We call it {\em the $j$-conjecture.} Putting $j=1$,  it asks if
 \begin{equation}\label{normConj}\big \|\Re((I-A)^{-1})\big \|   \geq \big \| (I+|A|)^{-1}\big \|.
\end{equation}

The results shown below   are  weaker than the
conjectured inequalities. 

\begin{pro} Let $A$ be an $n\times n$ strict contraction. Then,
for $j=1, 2, \dots, n$,
\begin{equation}\label{weak1}
 \l_j\big (\Re((I-A)^{-1})\big )\geq \l_j\big ((I+|A|)^{-1}\big )- \frac{r_{1}^2-r_{n-j+1}^2}{2(1+r_{n-j+1})^2}
\end{equation}
and
\begin{equation}\label{weak2}
 \big \|\Re((I-A)^{-1})\big \|   \geq \big \| (I+|A|)^{-1}\big \|-
\frac{r_{1}^2-r_{n}^2}{2(1+r_{n})^2}.
\end{equation}
 \end{pro}

   \proof We derive as follows.
 \begin{eqnarray*}
\l_j\big (\Re((I-A)^{-1})\big ) & = & \frac12
\l_j\big (I+(I-A^*)^{-1}(I-A^*A)(I-A)^{-1}\big ) \\
 & = & \frac12
\big (1+\l_j ( (I-A^*)^{-1}(I-A^*A)(I-A)^{-1})\big ) \\
& \geq  & \frac12
\big (1+\l_j ( (I-A)^{-1}(I-A^*)^{-1})\,\l_n(I-A^*A)\big ) \\
& =  & \frac12 \Big (1+
\frac{\l_n(I-A^*A)}{\l_{n-j+1}((I-A^*)(I-A))} \Big ) \\
&= & \frac12 \Big (1+
\frac{1-\l_{1}(A^*A)}{\s_{n-j+1}^2(I-A)} \Big ) \\
& \geq  & \frac12 \Big (1+
\frac{1-r_{1}^2}{(1+r_{n-j+1})^2} \Big ) \\
& = & \frac{1}{1+r_{n-j+1}}-  \frac{r_{1}^2-r_{n-j+1}^2}{2(1+r_{n-j+1})^2}\\
& = & \l_j\big ((I+|A|)^{-1}\big )- \frac{r_{1}^2-r_{n-j+1}^2}{2(1+r_{n-j+1})^2}.
\end{eqnarray*}

This completes the proof of (\ref{weak1}). Setting $j=1$ results in (\ref{weak2}). $\eop$

In a similar way, we can obtain,
for $j=1, 2, \dots, n$,
\begin{equation}
\l_j\big (\Re((I-A)^{-1})\big )\geq  \frac{1}{1+r_{1}}+  \frac{r_{1}^2-r_{n-j+1}^2}{2(1+r_{1})^2}
\end{equation}
and
\begin{equation}\label{eq20}
 \l_j\big ((I-A^*)^{-1}(I-A^*A)(I-A)^{-1}\big )    \geq \frac{1-r_{n-j+1}^2}{(1+r_{1})^2}.
\end{equation}

A few special cases of the open problem  have been settled.

(I). The $j$-conjecture 
holds true for normal contractions (including positive semidefinite matrices, Hermitian matrices), i.e., for $A$ with $\|A\|<1$ and $A^*A=AA^*$. This is due to the fact that  normal matrices are unitarily diagonalizable and that if $c$ is a complex number with $|c|<1$, then
$\Re ((1-c)^{-1})\geq (1+|c|)^{-1}.$

(II). The $j$-conjecture 
holds true for $j=n$
by (\ref{weak1}).

(III). The $j$-conjecture 
 holds true for $j=1$ and singular contractions, i.e., $r_n=0$.  Since $A$ is singular, there exists a unit vector $u$ such that $Au=0$. Observe that
$$(I-A)^{-1}=I+(I-A)^{-1}A.$$
We have
$$u^*\Re((I-A)^{-1}) u= \Re(u^*(I-A)^{-1}u) =1.$$
Since $\Re((I-A)^{-1})$ is positive definite, the spectral norm of $\Re((I-A)^{-1})$ is the same as its largest eigenvalue.
The min-max principle reveals at once
$$\|\Re((I-A)^{-1})\|=\l_1(\Re((I-A)^{-1}))=\max_{\|x\|=1}x^* \Re((I-A)^{-1}) x \geq 1=\| (I+|A|)^{-1}\|.$$

\begin{rem}
\rm With $(I-A)^{-1}=I+(I-A)^{-1}A$, we have
$$\Re((I-A)^{-1})=I+\Re((I-A)^{-1}A).$$
Thus, (\ref{normConj}) is equivalent to 
\begin{equation}\label{normConjRem}
\l_1\big (\Re((I-A)^{-1}A)\big )   \geq
-\frac{r_n}{1+r_n}.
\end{equation}
\end{rem}

\begin{rem} \rm  Note that $\|\Re(I-A)^{-1}\|=\max_{\|x\|=1} \Re  (x^*(I-A)^{-1}x )$.
For  $X\in M_n$,
since $\max_{\|x\|=1} x^*(\Re (X))x =\max_{\|x\|=1}\Re (x^*Xx)  \geq \max_{j} \Re (\l_j(X))$,
the norm inequality (\ref{normConj}) would follow from  the second inequality below
$$\max_{\|x\|=1} \Re \big (x^*((I-A)^{-1})x\big )\geq \max_j\Re \big  (\l_j ((I-A)^{-1})\big )\geq
\l_1\big ((I+|A|)^{-1}\big ).$$
That is, to show (\ref{normConj}), it   suffices  to prove that $A$ has an eigenvalue $\l$  such that
$$\Re \big ( (1-\l)^{-1}\big ) \geq (1+r_n)^{-1}.$$
However, this is not true in general.  Take
$$A =\left ( \begin{array}{rrr}
-0.2007  &  0.0263  &   -0.4910\\
    0.5055 &  -0.2419 &   0.5709\\
    0.3799  &  0.1640 &  -0.3848 \end{array} \right ).$$
    The eigenvalues of $A$ are
     $-0.1482 + 0.3451i,$ $
 -0.1482 - 0.3451i$,
  $-0.5309$, and the singular values of $A$ are
  $0.9554$, $
    0.5556$, $
    0.1411$.  Upon computation, we have
    $$\max_{\l\in \Gamma(A)} \Re\big ( (1-\l)^{-1}\big ) = 0.7988< 0.8763 = (1+r_n)^{-1}.$$
  Note that $\| \Re((I-A)^{-1})\|=1.0301.$
\end{rem}

 \section{Cayley transforms with majorization}

This section is devoted to the partial products of singular values of the Cayley transforms of   given matrices. Cayley transform is originally defined for real skew-symmetric matrices which have no nonzero real eigenvalues (as a result the Cayley transform matrix is orthogonal). To be precise, let $S$ be a real skew-symmetric matrix, then $\mathcal{C}(S)=(I+S)(I-S)^{-1}$ is called the
Cayley transform of $S$ (see, e.g.,  \cite[p.\,73]{GolubMC96} or \cite[p.\,75]{ZhanAMSBook13}).
For  linear operators on Hilbert spaces, there is a rich theory about Cayley transform  with
 linear dissipative operator, contraction, and isometry.
Let $X\in M_n$. If $X+iI$ is invertible,
 we call
$\mathcal{C}(X)=(X-iI)(X+iI)^{-1}$ the {\em Cayley transform} of $X$.
(More generally,  for a nonsingular matrix $A$,  $A^{-1}A^*$ is called {\em generalized Cayley transform} of $A$ \cite{FanLAA1972}.)
A large family of matrices with well-defined Cayley transforms  exists: strict contractions, positive semidefinite matrices, Hermitian matrices,
stable matrices, and matrices with all real eigenvalues, etc. We are concerned with the
Cayley transforms of contractions.

\begin{thm}\label{Cay}
Let $A, B\in M_n$ be strict contractions and let $\mathcal{C}(A)$ and $\mathcal{C}(B)$ be the Cayley transforms of $A$ and $B$, respectively. Then for $1\leq i_1<\cdots <i_k\leq n$,
{\small $$\prod_{j=1}^k \frac{1-\s_{n-i_j+1}(A)}{1+\s_{j} (A)}\leq \prod_{j=1}^k \s_{i_j}\big (\mathcal{C}(A)\big )\leq  \prod_{j=1}^k \frac{1+\s_{i_j}(A)}{1-\s_{j}(A)}$$}
and
{\small $$\prod_{j=1}^k \frac{2\,\s_{i_j}(A-B)}{(1-\s_{j} (A))(1-\s_j(B))}\leq \prod_{j=1}^k \s_{i_j}\big (\mathcal{C}(A)-\mathcal{C}(B)\big ) \leq \prod_{j=1}^k \frac{2\,\s_{i_j}(A-B)}{(1-\s_{j} (A))(1-\s_j(B))}.$$}
\end{thm}

\proof  We compute the upper bounds. The lower bounds are similarly derived.
\begin{eqnarray*}
\prod_{j=1}^k \s_{i_j}\big (\mathcal{C}(A)\big ) & = & \prod_{j=1}^k \s_{i_j}\big ((A-iI)(A+iI)^{-1} \big )\\
& \leq & \prod_{j=1}^k \s_{i_j} (A-iI)\, \s_j\big ((A+iI)^{-1} \big )\\
& \leq & \prod_{j=1}^k \s_{i_j} (A-iI) \big (\s_{n-j+1} (A+iI) \big )^{-1}\\
& \leq & \prod_{j=1}^k \frac{1+\s_{i_j}(A)}{1-\s_{j}(A)}.
\end{eqnarray*}

We used fact (iv) in the above derivation. The lower bound is obtained by using (iii).
For the second part, we observe  that $\mathcal{C}(A)=I-2i(A+iI)^{-1}$ and
$$\mathcal{C}(A)-\mathcal{C}(B)= 2i(B+iI)^{-1}(A-B)(A+iI)^{-1}.$$
It follows that, by using (iv) twice,
\begin{eqnarray*}
\prod_{j=1}^k \s_{i_j}\big (\mathcal{C}(A)-\mathcal{C}(B)\big ) & = &  \prod_{j=1}^k 2\,\s_{i_j}\big ((B+iI)^{-1}(A-B)(A+iI)^{-1} \big )\\
& \leq & \prod_{j=1}^k 2\, \s_{j}\big ((A+iI)^{-1}\big)\, \s_{i_j}(A-B)\, \s_j\big ((B+iI)^{-1}\big ) \\
& =& \prod_{j=1}^k \frac{2 \,\s_{i_j}(A-B)}{\s_{n-j+1}(A+iI)\, \s_{n-j+1} (B+iI) }\\
& \leq & \prod_{j=1}^k \frac{2\,\s_{i_j}(A-B)}{(1-\s_{j}(A))\, (1-\s_j(B))}. \eop
\end{eqnarray*}

 Setting $k=1$ in the theorem, we obtain the lower and upper bounds for the singular values of the  Cayley transforms of strict contractions $A$, that is,
$$  \frac{1-\s_{n-j+1}(A)}{1+\s_{1} (A)}\leq   \s_{j}\big (\mathcal{C}(A)\big )\leq  \frac{1+\s_{j}(A)}{1-\s_{1}(A)}, \;\; j=1, 2, \dots, n.$$

\begin{rem} \rm
The proof of Theorem \ref{Cay} was in the spirit of Fan and Hoffman's \cite{FanHoffman55} in which stronger inequalities were shown with $A$ and $B$ being 
 Hermitian: $
\s_{j}\big (\mathcal{C}(A)-\mathcal{C}(B)\big )\leq  2\,\s_{j}(A-B).$
We  point out that
a weaker version of the Fan and Hoffman result is stated in \cite[p.\,374]{MOA11}
as $\|A-B\|\geq \frac12 \|\mathcal{C}(A)-\mathcal{C}(B)\|$ for all unitarily invariant norms which is equivalent to
the weak majorization
$\s(A-B) \succ_w  \frac12  \s(\mathcal{C}(A)-\mathcal{C}(B))$. (Note: there is a typo in the display (12b) in the book, i.e.,
$\prec_w$ should be $\succ_w$.) Our results are given as log-majorization (which implies weak majorization; see, e.g., \cite[p.\,345]{ZFZbook11}) for more general matrices.
\end{rem}

\section{Multiple matrices}

Let $A, B\in M_n$ be such that $A$ and $A-B$ are nonsingular. One can check that
$$(A^*-B^*)^{-1}(A^*A-B^*B)(A-B)^{-1}=2\, \Re ((I-BA^{-1})^{-1})-I.$$

If, additionally,  the spectral norm $\|BA^{-1}\|<1$, then,  by (i), for each $j$,
\begin{eqnarray*}
\lefteqn{\l_{j} \big ((A^*-B^*)^{-1}(A^*A-B^*B)(A-B)^{-1}\big )} \\
& & = \l_{j} \big ( 2\, \Re ( (I-BA^{-1})^{-1} ) -I\big ) \\
& & \leq  2\,\s_{j} (  (I-BA^{-1})^{-1} ) -1   \\
& & = \frac{2}{\s_{n-{j}+1} (  I-BA^{-1})} -1 \\
& & \leq  \frac{2}{1-\s_{j}(BA^{-1} )} -1 \\
& & \leq  \frac{2}{1-\s_{j} (B)\s_1(A^{-1} )} -1 \\
&  & \leq  \frac{\s_n(A)+\s_j(B)}{\s_n(A)-\s_{j}(B)}.
\end{eqnarray*}  
Consequently, for
contractions $A$ and $B$ with $\|BA^{-1}\|<1$ and   $\det(A-B)\not =0$,
$$
 \prod_{j=1}^k\l_{i_j} \Big ((A^*-B^*)^{-1}(A^*A-B^*B)(A-B)^{-1}\Big ) \leq
 \prod_{j=1}^k\frac{\s_n(A)+\s_j(B)}{\s_n(A)-\s_{j}(B)}.$$


\subsection*{Acknowledgments} {\small   This work
 was done when Chaojun Yang was a CSC-sponsored Ph.D. student at Nova Southeastern University during the  2019-2020 academic year. The authors thank Prof.~Zhaolin Jiang for initiating the work and Prof.~Lei Cao for discussions.}

\end{document}